\input amstex
\documentstyle{amsppt}\nologo\footline={}

\def\Isom{\mathop{\text{\rm Isom}}}
\def\Cl{\mathop{\text{\rm Cl}}}
\def\T{\mathop{\text{\rm T}}}
\def\N{\mathop{\text{\rm N}}}

\input epsf
\input lpic-mod

\hsize450pt\topmatter\title Yet Another Poincar\'e's Polyhedron
Theorem\endtitle\author Sasha Anan$'$in and Carlos
H.~Grossi\endauthor\address Departamento de Matem\'atica, IMECC,
Universidade Estadual de Campinas,\newline13083-970--Campinas--SP,
Brasil\endaddress\email Ananin$_-$Sasha\@yahoo.com\endemail\address
Max-Planck-Institut f\"ur Mathematik, Vivatsgasse 7, 53111 Bonn,
Germany\endaddress\email
grossi$_-$ferreira\@yahoo.com\endemail\subjclass 51M10 (22E40,
51H05)\endsubjclass\abstract Poincar\'e's Polyhedron Theorem is a
widely known valuable tool in constructing manifolds endowed with a
prescribed geometric structure. It is one of the few criteria providing
discreteness of groups of isometries. This work contains a version of
Poincar\'e's Polyhedron Theorem that is applicable to constructing
fibre bundles over surfaces and also suits geometries of nonconstant
curvature. Most conditions of the theorem, being as local as possible,
are easy to verify in practice.\endabstract\endtopmatter

\document

\centerline{\bf1.~Introduction}

\medskip

{\bf1.1.}~It is frequently important to decide if a given subgroup $G$
of a Lie group is discrete. For instance, such a necessity appears
while constructing manifolds endowed with a prescribed geometric
structure. Typically, the group $G$ is related to some geometrical
configuration: it acts on a simply-connected homogeneous space $M$ and
is generated by isometries that identify given codimension $1$
subspaces called faces. These faces may bound a (fundamental)
polyhedron $P$ in $M$ into which the quotient $M/G$ can be `cut and
unfolded.' Thus, we expect certain pairs of faces to be identified by
the generators of $G$, called face-pairing isometries, in such a way
that $M/G$ results from the identifications and the space $M$ is
tessellated by the copies of $P$. This can be reversed: starting with a
polyhedron that tessellates $M$, we get the discrete group generated by
the face-pairing isometries. We have just briefly described the general
settings surrounding Poincar\'e's Polyhedron Theorem (PPT). The theorem
has a long history and plenty of versions; the interested reader may
consult, for instance, [EPe] and the references therein.

The main step in verifying the tessellation of $M$ is usually the study
of tessellation around the codimension $2$ faces of $P$, called edges.
This leads to the concept of a geometric cycle of edges: given an edge
$e$, its geometric cycle is a cyclic sequence of edges related by
face-pairing isometries such that the corresponding copies of $P$ are
expected to tessellate $M$ around $e$. Dealing with geodesic polygons
in the hyperbolic plane, Henri Poincar\'e realized that in order to
obtain the tessellation of $M$ it suffices to require that the sum of
the angles of the polygon along every geometric cycle equals $2\pi$
(for simplicity, we do not deal here with ideal cycles). Later, he
extended this idea to the case of constant curvature hyperbolic
$3$-space.

\smallskip

{\bf1.2.}~Most versions of PPT come from constant curvature (or even
plane) geometries, where convexity arguments play an important role
typically suited to polyhedra with constant angles between (totally
geodesic) faces along common edges. In general, such an approach is
inapplicable to nonconstant curvature geometries, say, to complex
hyperbolic geometry. The usual requirements like `adjacent polyhedra
intersect in an expected way' are difficult to check. Therefore, we
look for a version where the conditions for tessellation are as local
as possible and provide global properties just {\it a posteriori.} The
strategy is to impose some requirements of infinitesimal nature which
can be verified in practice and then obtain an infinitesimal
tessellation that can be `integrated' with the help of suitable local
conditions expressing a good behaviour of the faces. Note that while
constructing manifolds we are used to having an explicitly given set of
face-pairing isometries. We therefore treat the relation between the
face-pairing isometries involved in a cycle of edges as being easily
verifiable: at worst, we need to multiply a few matrices.

\smallskip

{\bf1.3.}~We show that the tessellation of a metric neighbourhood (see
Tessellation Condition 2.1) is~sufficient for discreteness (Proposition
2.2). This condition seems to be more useful than the well-known
completeness requirement. We think it can be particularly relevant in
dimension $>2$ if one also has to deal with the `parabolic cycles.'
Conditions similar to Tessellation Condition 2.1 have already appeared
in the literature (see, for instance, [Ale] and [Bea]).

The Poincar\'e angle condition can be weakened to the following total
angle condition. Pick a point $p$ in an edge $e$. Being subsequently
applied to $p$, the face-pairing isometries involved in the geometric
cycle of $e$ provide a point in each edge in the cycle. We measure the
interior angles between faces at these points. The total angle at $p$,
i.e., the sum of such angles, is a multiple of $2\pi$ because of the
cycle relation. {\it Total angle condition\/} means that, for every
cycle, the total angle at a single point equals $2\pi$. Such a
condition can be quite handy in particular cases: sometimes the
verification happens to be very simple at geometrically distinguished
points.

In Theorem 3.5, we deal with polyhedra that possess no faces of
codimension $>2$. In this case, total angle condition essentially
ensures the tessellation of a topological neighbourhood of the
polyhedron. Requiring in addition the metric separability of faces
(including that of faces sharing a common edge; see~Strong Simplicity
IV and Condition (3) in Theorem 3.5) allows us to integrate the
topological tessellation into the tessellation of a metric
neighbourhood of the polyhedron. In Final Remarks 4.3, we explain how
one can in principle use the ideas of [Ale] and the current paper to
obtain a more general form of the theorem with no restriction on
codimension of faces. As it stands, Theorem 3.5 is well adapted to the
construction of fibre bundles over surfaces (in order of importance,
such manifolds are probably the first after the compact ones); for
example, it applies directly to constructing complex hyperbolic disc
bundles in [AGG].

A serious defect of our version is the global requirement of
simplicity, i.e., the requirement that the faces intersect as expected
and thus bound the polyhedron itself. In complex hyperbolic geometry,
for example, it is already difficult to check the simplicity of a
polyhedron with bisectors taken as faces. Hence, it seems that one
should obtain an even more local PPT which makes the verification of
simplicity unnecessary. This would finally `disassemble' the
polyhedron, taking away the arbitrary choice involved in PPT.

\medskip

{\bf Acknowledgments.} We are grateful to Pedro Frejlich for his
constant support, to Misha Kapovich for discussions on the subject, and
to the referee for valuable suggestions. We are indebted to the IH\'ES
for hospitality, support, and a research style that made it possible to
accomplish many of our scientific plans, including this work.

\bigskip

\centerline{\bf2.~Preliminaries}

\medskip

This is essentially standard material (see, for instance, [Bea, \S9.8,
p.~242]).

\smallskip

Let $M$ be a locally path-connected, connected, and simply-connected
metric space. Denote by $B(x,\varepsilon)$ the open ball of radius
$\varepsilon>0$ centred at $x$ and let
$N(X,\varepsilon):=\bigcup\limits_{x\in X}B(x,\varepsilon)$ for
$X\subset M$. We regard a {\it polyhedron\/} in $M$ as being a closed,
locally path-connected, and connected subspace $P\subset M$ such that

\smallskip

$\bullet$ $P$ is the closure of its nonempty interior:
$\overset{\,_\circ}\to P\ne\varnothing$ and
$P=\Cl\overset{\,_\circ}\to P$;

\smallskip

$\bullet$ the nonempty boundary of $P$ is decomposed into the union of
nonempty subsets $s\in S$ called {\it faces\/}:
$\partial P:=P\setminus\overset{\,_\circ}\to P=\bigcup\limits_{s\in
S}s$.

\smallskip

A {\it face-pairing\/} of a polyhedron $P$ is an involution
$\overline{\phantom{m}}:S\to S$ and a family of isometries
$I_s\in\Isom M$ satisfying $I_ss=\overline s$ and
$I_{\overline s}=I_s^{-1}$ for every face $s\in S$.

Let $P$ be a polyhedron with a given face-pairing and let $G$ denote
the group generated by the face-pairing isometries. We introduce a
relation in $G\times P$ by putting $(g,x)\sim(h,y)$ exactly when
$x\in s$ for some $s\in S$, $I_sx=y$, and $h^{-1}g=I_s$. Taking the
closure of this symmetric relation with respect to transitivity (and
reflexivity), we obtain an equivalence relation also denoted by $\sim$.
Let $J:=G\times P/\sim$ and let $[g,x]$ denote the class of $(g,x)$ in
$J$. Consider the discrete topology on $G$ and equip $P$, $G\times P$,

{\unitlength=1bp$$\latexpic(0,0)(-170,38)
\put(0,40){$G\times P$}\put(30,44){\vector(1,0){15}}\put(34,47){$\pi$}
\put(47,40){$J$}\put(13,35){\vector(1,-2){12}}\put(8,20){$\psi$}
\put(45,35){\vector(-1,-2){12}}\put(43,20){$\varphi$}\put(24,0){$M$}
\endlatexpic$$}

\rightskip67pt

\vskip-43pt

\noindent
and $J$ with their natural topologies. We have a commutative diagram of
continuous $G$-maps $\psi(g,x):=gx$, $\pi(g,x):=[g,x]$, and
$\varphi[g,x]:=gx$. (Actions of $G$ by homeomorphisms are defined by
$h(g,x):=(hg,x)$ and $h[g,x]:=[hg,x]$.) Let

\vskip12pt

\centerline{$[P]:=\big\{[1,x]\mid x\in P\big\}\quad\text{\rm
and}\quad[\overset{\,_\circ}\to P]:=\big\{[1,x]\mid
x\in\overset{\,_\circ}\to P\big\}$.}

\vskip10pt

\rightskip0pt

\noindent
Clearly, $J=\bigcup\limits_{g\in G}g[P]$ and
$g_1[\overset{\,_\circ}\to P]\cap g_2[\overset{\,_\circ}\to
P]\ne\varnothing$
implies $g_1=g_2$. In other words, $[P]$ is a {\it fundamental
region\/} for the action of $G$ on $J$.

We assume that $\pi^{-1}[1,x]$ is {\bf finite} for every
$x\in\partial P$, hence, for every $x\in P$. Let $x\in P$. Then
$\pi^{-1}[1,x]=\big\{(g_1,x_1),\dots,(g_n,x_n)\big\}$ for some
$g_j\in G$ and $x_j\in P$. The polyhedra $g_jP$ are the {\it formal
neighbours\/} of $P$ at $x$. For $\delta>0$, define
$$N_{x_j,\delta}:=\big\{y\in P\mid d(y,x_j)<\delta\big\}\subset
P,\qquad N_{x,\delta}:=\bigcup_{j=1}^n(g_j,N_{x_j,\delta})\subset
G\times P,\qquad W_{x,\delta}:=\pi N_{x,\delta}\subset J,$$
where $d(\cdot,\cdot)$ stands for the distance function on $M$. Using
this notation, we state the

\medskip

{\bf2.1.~Tessellation Condition.} A polyhedron $P$ with a given
face-pairing satisfies Tessellation Condition if

\smallskip

$\bullet$ for every $x\in P$, there exists some $\delta(x)>0$ such that
$\pi^{-1}(W_{x,\delta})=N_{x,\delta}$ and
$\varphi W_{x,\delta}=B(x,\delta)$ for all $0<\delta\le\delta(x)$;

$\bullet$ some open metric neighbourhood $N$ of $P$ in $M$ is
tessellated; this means that $N(P,\varepsilon)\subset N$ for some
$\varepsilon>0$ and that there exists a function $f:P\to\Bbb R$ taking
positive values such that $\varphi:W_{P,f}\to N$ is bijective, where
$W_{P,f}:=\bigcup\limits_{x\in P}W_{x,f(x)}$.

\medskip

{\bf2.2.~Proposition.} {\sl Tessellation Condition\/ {\rm2.1} implies
that\/ $\varphi$ is a homeomorphism. In other words, the polyhedron\/
$P$ is a fundamental region for the action of\/ $G$ on\/ $M$.}

\medskip

{\bf Proof.} Straightforward arguments show that $J$ is Hausdorff and
path-connected, that the family
$\big\{gW_{x,\delta}\mid g\in G,\,x\in P,\,0<\delta\le\delta(x)\big\}$
is a base of the topology on $J$, and that $\varphi$ is a local
homeomorphism. Clearly, $\varphi:gW_{P,f}\to gN$ is a homeomorphism for
all $g\in G$. As $M$ is simply-connected, it~suffices to show that
$\varphi$ is a regular covering.

Since $\varphi$ is open, $\varphi J$ is open in $M$. Let
$x\in\Cl(\varphi J)$. Then $B(x,\varepsilon)\cap gP\ne\varnothing$ for
some $g\in G$. It~follows that
$x\in N(gP,\varepsilon)\subset gN=\varphi(gW_{P,f})\subset\varphi J$.
Hence, $\varphi J$ is closed in $M$. Since $M$ is connected, $\varphi$
is surjective.

Take $x\in M$. Define
$$G_x:=\big\{g\in G\mid U_x\cap gP\ne\varnothing\big\},$$
where $U_x\subset B(x,\frac12\varepsilon)$ is a path-connected open
neighbourhood of $x$. For every $g\in G_x$, let
$$W_g:=\varphi^{-1}(U_x)\cap gW_{P,f}.$$
Since $U_x\cap gP\ne\varnothing$ implies that
$U_x\subset B(x,\frac12\varepsilon)\subset N(gP,\varepsilon)\subset
gN$,
we conclude that $\varphi:W_g\to U_x$ is a homeomorphism. Moreover,
$$\varphi^{-1}(U_x)=\bigcup_{g\in G_x}W_g.$$
It remains to show that the distinct $W_g$'s are disjoint. Suppose that
$W_{g_1}\cap W_{g_2}\ne\varnothing$ for some $g_1,g_2\in G_x$. The
projection $W_{g_1}\times W_{g_2}\to W_{g_1}$ induces a homeomorphism
between
$$X:=\big\{(x_1,x_2)\in W_{g_1}\times W_{g_2}\mid\varphi x_1=\varphi
x_2\big\}$$
and $W_{g_1}$. The diagonal
$$\Delta_{W_{g_1}\cap W_{g_2}}=\Delta_J\cap(W_{g_1}\times
W_{g_2})\subset X$$
is closed in $X$ since $J$ is Hausdorff. Therefore, the image
$W_{g_1}\cap W_{g_2}$ of $\Delta_{W_{g_1}\cap W_{g_2}}$ is closed in
$W_{g_1}$. Since $W_{g_1}$ is connected, we obtain $W_{g_1}=W_{g_2}$
$_\blacksquare$

\bigskip

\centerline{\bf3.~A Plane-like Poincar\'e's Polyhedron Theorem}

\medskip

In what follows, $M$ is a connected, oriented, and simply-connected
Riemannian manifold. We regard a {\it cornerless polyhedron\/}
$P\subset M$ with a face-pairing as a subspace satisfying the
conditions stated in the beginning of the previous section as well as
those below.

\medskip

{\bf I.}~The faces of $P$ are topologically closed, oriented smooth
connected submanifolds of codimension~$1$ in $M$ with (possibly empty)
boundary. Each face $s$ of $P$ is oriented so that normal vectors to
$s\setminus\partial s$ point towards the interior of $P$.

\smallskip

{\bf II.}~The boundary of every face $s\in S$ is a disjoint union
$\partial s=\bigsqcup\limits_{e\in E_s}e$ of nonempty connected {\it
edges.} ($E_s=\varnothing$ is allowed.) We write $e\diamond s$ or
$s\diamond e$ if $e\in E_s$. Clearly, $e\diamond s$ implies
$\overline s\diamond I_se$.

\smallskip

{\bf III.}~$P$ has a finite number of faces and edges. Each edge $e$
belongs to exactly two distinct faces $s_1$ and $s_2$. In symbols:
$s_1\diamond e\diamond s_2$.

\smallskip

{\bf IV~{\rm(Strong Simplicity)}.} The intersection of two distinct
faces is contained in the boundary of both faces and is a (possibly
empty) union of edges. The distances between:

\quad two distinct edges,

\quad two distinct faces that do not share an edge,

\quad a face and an edge not contained in it

\noindent
are all greater than some $d>0$.

\medskip

{\bf3.1.}~Start with $\overline s_0\diamond e\diamond s_1$. Applying
$I_{s_1}$ to $s_1$ and $e$, we obtain
$\overline s_1\diamond I_{s_1}e\diamond s_2$. Applying $I_{s_2}$ to
$s_2$ and~$I_{s_1}e$, we obtain
$\overline s_2\diamond I_{s_2}I_{s_1}e\diamond s_3$, and so on. (Of
course, by III, we eventually arrive back at
$\overline s_0\diamond e\diamond s_1$.)

A cyclic sequence
$$\overset{^{I_{s_n}}}\to\longrightarrow\overline s_n=\overline
s_0\diamond e\diamond
s_1\overset{^{I_{s_1}}}\to\longrightarrow\overline s_1\diamond
I_{s_1}e\diamond s_2\overset{^{I_{s_2}}}\to\longrightarrow\overline
s_2\diamond I_{s_2}I_{s_1}e\diamond
s_3\overset{^{I_{s_3}}}\to\longrightarrow\dots\overset{^{I_{s_{n-1}}}}
\to\longrightarrow\overline s_{n-1}\diamond I_{s_{n-1}}\cdots
I_{s_1}e\diamond s_n\overset{^{I_{s_n}}}\to\longrightarrow,$$
where each term is obtained from the previous one by the above rule, is
called a {\it cycle of edges.} The~number $n$ is the length of the
cycle and the isometry $I:=I_{s_n}\cdots I_{s_1}$ will be referred to
as the {\it cycle isometry.} A cycle can be read backwards, i.e., in
opposite {\it orientation,} which inverts its isometry. If the cycle
isometry is the identity and if the cycle is the shortest one with this
property, then the cycle is said to be {\it geometric\/} (see also
Remark 3.4.) Clearly, every cycle is a multiple of a shortest,
combinatorial one. Note that, in a geometric cycle, an edge may occur
several times (this does not happen in a combinatorial cycle).

Assume that we are given a family of disjoint geometric cycles that
contains every edge of $P$. Fixing a term (say)
$\overline s_0\diamond e\diamond s_1$ in some oriented cycle of the
family, define $I_j:=I_{s_j}\cdots I_{s_1}$ for all $j=0,1,\dots,n$ (we
usually consider $j$ modulo $n$) so that the cycle takes the form
$$\overset{^{I_{s_n}}}\to\longrightarrow\overline s_n=\overline
s_0\diamond I_0e\diamond
s_1\overset{^{I_{s_1}}}\to\longrightarrow\overline s_1\diamond
I_1e\diamond s_2\overset{^{I_{s_2}}}\to\longrightarrow\overline
s_2\diamond I_2e\diamond s_3\overset{^{I_{s_3}}}\to\longrightarrow\dots
\overset{^{I_{s_{n-1}}}}\to\longrightarrow\overline s_{n-1}\diamond
I_{n-1}e\diamond s_n\overset{^{I_{s_n}}}\to\longrightarrow.$$

{\bf3.2.}~We can describe all formal neighbours of $P$ at a point
$x\in\partial P$. If $x$ does not belong to any edge, then there is a
unique face $s$ containing $x$,
$\pi^{-1}[1,x]=\big\{(1,x),(I_{\overline s},I_sx)\big\}$, and the only
formal neighbours of $P$ at $x$ are $P$ and $I_{\overline s}P$. If $x$
belongs to an edge $\overline s_0\diamond e\diamond s_1$, then
$\pi^{-1}[1,x]=\big\{(I_j^{-1},I_jx)\mid
j=0,1,\dots,n-\nomathbreak1\big\}$
and the $I_j^{-1}P$ are the formal neighbours of $P$ at $x$. Indeed,
suppose that $(I_j^{-1},I_jx)\sim(h,y)$. This means that $I_jx\in s'$,
$I_{s'}I_jx=y$, and $h^{-1}I_j^{-1}=I_{s'}$ for some $s'\in S$. In
particular, $I_je$~and $s'$ intersect. It follows from IV that
$I_je\diamond s'$. Hence, either $s'=\overline s_j$ or $s'=s_{j+1}$.
Therefore, either $(h,y)=(I_{j-1}^{-1},I_{j-1}x)$ or
$(h,y)=(I_{j+1}^{-1},I_{j+1}x)$. It remains to observe that the $I_j$,
$j=0,1,\dots,n-1$, are all distinct because we could otherwise take a
shorter cycle whose isometry would be the identity.

\smallskip

{\bf3.3.}~Pick a point $x$ in some edge
$\overline s_0\diamond e\diamond s_1$ of an oriented cycle. Let
$\N_xe:=(\T_xe)^\perp$ and let $n_0,n_1$ denote, respectively, the unit
normal vectors to $\overline s_0,s_1$ at $x$ that point towards the
interior of $P$. Let~$t_0\in\T_x\overline s_0\cap\N_xe$ and
$t_1\in\T_xs_1\cap\N_xe$ stand for the unit vectors that point
respectively towards the interiors of $\overline s_0$ and $s_1$. The
basis $t_0,n_0$ orients $\N_xe$. This orientation corresponds to the
orientation of the\break

\vskip-5pt

\noindent
\hskip341pt$\vcenter{\hbox{\epsfbox{Picture1.eps}}}$

\rightskip120pt

\vskip-55pt

\noindent
cycle. The oriented {\it interior angle\/} $\alpha_0$ from
$\overline s_0$ to $s_1$ at $x$ is the angle from $t_0$ to $t_1$, which
takes values in $[0,2\pi]$. We define similarly the interior angle
$\alpha_j$ from $\overline s_j$ to $s_{j+1}$ at $I_jx$. The sum
$\sum_{j=0}^{n-1}\alpha_j$ is the {\it total interior angle\/} of the
cycle at $x$. It is easy to see that altering the orientation of the
cycle alters the orientation of the corresponding $\N_xe$ and keeps the
same values of the~$\alpha_j$'s.

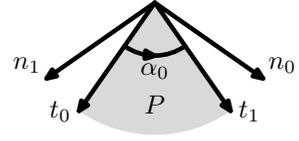
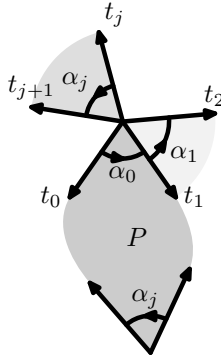
\rightskip0pt

\vskip7pt

\noindent
$\vcenter{\hbox{\epsfbox{Picture2.eps}}}$

\leftskip90pt

\vskip-141pt

Suppose that the face-pairing isometries {\it send interior into
exterior.} By definition, this means that $I_sn_s=-n_{\overline s}$ for
every face $s\in S$, where $n_s$ stands for the unit normal vector to
$s$ at some $x\in s$. This property implies the following. Take a point
$x$ in some edge $\overline s_0\diamond e\diamond s_1$ of an oriented
geometric cycle. Let~$t_1\in\T_xs_1\cap\N_xe$ be the unit vector that
points towards the interior of $I_1^{-1}\overline s_1=s_1$ and let
$t_2\in\T_xI_1^{-1}s_2\cap\N_xe$ be the unit vector that points towards
the interior of~$I_1^{-1}s_2$. Then the oriented angle from $t_1$ to
$t_2$ equals $\alpha_1$. In the same way, denoting by
$t_j\in\T_xI_j^{-1}\overline s_j\cap\N_xe$ the unit vector that points
towards the interior of~$I_j^{-1}\overline s_j$, we~can see that the
oriented angle from $t_j$ to $t_{j+1}$ equals $\alpha_j$. This implies
immediately that $\sum_{j=0}^{n-1}\alpha_j\equiv0\mod2\pi$. In
particular, the total interior angle of a geometric cycle is {\bf
constant\/}: it does not depend on the choice of $x\in e$.

Obviously, the distinct formal neighbours of $P$ at a point in an edge
overlap\penalty-10000

\leftskip0pt

\vskip-12pt

\noindent
when the total interior angle of a cycle is different from $2\pi$. In
the terms of Proposition 2.2, this corresponds to a ramification of
$\varphi$.

\medskip

{\bf3.4.~Remark.} For some geometries, the nature of edges allows to
(formally) weaken the condition that the cycle isometry is the
identity. This happens in the case when every isometry $I$ that fixes
pointwise some edge $e$ is completely determined by the rotation angle
about some $x\in e$, that is, by the image $In\in\N_xe$ of some
$0\ne n\in\N_xe$. In this case, it suffices to require only that
$I|_e=1_e$ and that the total interior angle at $x$ vanishes modulo
$2\pi$.

\medskip

{\bf3.5.~Theorem.} {\sl Let\/ $P$ be a cornerless polyhedron with a
face-pairing providing a family of geometric cycles that contains every
edge of\/ $P$. Suppose that

\smallskip

{\bf(1)} the face-pairing isometries send interior into
exterior\/{\rm;}

{\bf(2)} the total interior angle equals\/ $2\pi$ at some point of an
edge for every cycle of the family\/{\rm;}

{\bf(3)} for every two distinct faces\/ $s,s'$ such that\/
$s\cap s'\ne\varnothing$ and for every\/ $\vartheta>0$, there exists\/
$\varepsilon=\varepsilon(s,s',\vartheta)>0$ such that\/
$s'\cap N(s,\varepsilon)\subset\bigcup\limits_{s\diamond e\diamond
s'}N(e,\vartheta)$.

\smallskip

\noindent
Then Tessellation Condition\/ {\rm2.1} is satisfied.}

\medskip

{\bf Proof.} In what follows, we denote $\widetilde X:=X\cap P$ for
$X\subset M$.

\medskip

{\bf First step.} Using Conditions (1--2), we will integrate (employing
IV) an infinitesimal tessellation into a topological one. So, we will
show that there exists a sufficiently small tessellated open ball
centred at $x$ for every $x\in P$. In other words, for every $x\in P$,
we will find some $\delta(x)>0$ such that the first part of
Tessellation Condition 2.1 is valid and, additionally,
$\varphi:W_{x,\delta}\to B(x,\delta)$ is injective for all
$0<\delta\le\delta(x)$. We distinguish the cases
$x\in s\setminus\partial s$ for some $s\in S$ and $x\in e$ for some
edge $e$.

\smallskip

$\bullet$ In the first case, choose $\delta_1>0$ such that
$B(x,\delta_1)$ does not intersect the edges of $s$ and such that
$B(x,\delta_1)\cap\partial P=B(x,\delta_1)\cap s$. Choose $\delta_2>0$
analogously with respect to $I_sx\in\overline s$. Let
$$N_{x,\delta}:=\big(1,\widetilde
B(x,\delta)\big)\bigcup\big(I_{\overline s},\widetilde
B(I_sx,\delta)\big),\qquad W_{x,\delta}:=\pi N_{x,\delta},$$
where $0<\delta\le\delta(x):=\min(\delta_1,\delta_2)$. Clearly,
$\pi^{-1}W_{x,\delta}=N_{x,\delta}$. We need to show that
$\varphi:W_{x,\delta}\to B(x,\delta)$ is a bijection.

Note that
$s\cap B(x,\delta)\subset\widetilde B(x,\delta)\cap I_{\overline
s}\widetilde B(I_sx,\delta)$.
Also,
$\widetilde B(x,\delta)\ne I_{\overline s}\widetilde B(I_sx,\delta)$ by
Condition (1). Pick a point
$q_0\in\widetilde B(x,\delta)\setminus I_{\overline s}\widetilde
B(I_sx,\delta)$
such that $q_0\notin s$. Due to the fact that $\delta\le\delta(x)$, a
smooth oriented curve $\gamma\subset B(x,\delta)$ connecting $q_0$ and
$q\in B(x,\delta)\setminus s$ can intersect $\partial P$ and
$\partial I_{\overline s}P$ only along
$(s\setminus\partial s)\cap B(x,\delta)$. We can assume that such
intersections are transverse. According to (1), when intersecting $s$,
the curve $\gamma$ leaves $\widetilde B(x,\delta)$ and enters
$I_{\overline s}\widetilde B(I_sx,\delta)$ or {\it vice-versa.} Hence,
$q$ belongs to exactly one of $\widetilde B(x,\delta)$ and
$I_{\overline s}\widetilde B(I_sx,\delta)$. The result then follows.

\smallskip

$\bullet$ The second case is similar. Let the $I_j$'s be related to the
geometric cycle containing $e$. Choose $\delta_j>0$ such that
$B(I_jx,\delta_j)$ does not intersect any edge of $\overline s_j$ or
$s_{j+1}$ except $I_je$ and such that
$B(I_jx,\delta_j)\cap\partial P=\big(B(I_jx,\delta_j)\cap\overline
s_j\big)\cup\big(B(I_jx,\delta_j)\cap s_{j+1}\big)$.
Let
$$N_{x,\delta}:=\bigcup_j\big(I_j^{-1},\widetilde
B(I_jx,\delta)\big),\qquad W_{x,\delta}:=\pi N_{x,\delta},$$
where $0<\delta\le\delta(x):=\min\delta_j$. The description 3.2 of
formal neighbours implies that $\pi^{-1}(W_{x,\delta})=N_{x,\delta}$.
We have
$$I_j^{-1}s_{j+1}\cap B(x,\delta)\subset I_j^{-1}\widetilde
B(I_jx,\delta)\cap I_{j+1}^{-1}\widetilde B(I_{j+1}x,\delta).$$
Let $F:=\bigcup\limits_jI_j^{-1}s_{j+1}\cap B(x,\delta)$ and let
$q_0\in\widetilde B(x,\delta)\setminus F$. Since $\delta\le\delta(x)$,
a smooth oriented $\gamma\subset B(x,\delta)$ connecting $q_0$ and
$q\in B(x,\delta)\setminus F$ may intersect
$\bigcup\limits_j\partial I_j^{-1}P$ only along $F$. We assume that
$\gamma$ does not intersect $e$ and is transverse to $F$. Condition (1)
implies that, when intersecting $I_j^{-1}s_{j+1}$, the curve $\gamma$
leaves $I_j^{-1}\widetilde B(I_jx,\delta)$ and enters
$I_{j+1}^{-1}\widetilde B(I_{j+1}x,\delta)$ or {\it vice-versa.} Hence,
$\varphi:W_{x,\delta}\to B(x,\delta)$ is surjective. Following the
discussion 3.3 concerning the total angle of the cycle at $x$, we
consider the closed sectors $T_j\subset\N_xe$ containing the oriented
interior angle of $I_j^{-1}P$ at $x$. Conditions (1) and (2) imply that
$\bigcup\limits_jT_j=N_xe$ and
$\overset{\circ}\to T_{j_1}\cap\overset{\circ}\to T_{j_2}=\varnothing$
if $j_1\not\equiv j_2\mod n$. Hence, distinct formal neighbours
$I_j^{-1}P$ cannot be equal.

Suppose that $\varphi:W_{x,\delta}\to B(x,\delta)$ is not injective at
some $q\in B(x,\delta)$. It follows from the description~3.2 of formal
neighbours that $q\notin F$. Pick a point $q_0$ living in exactly one
of the $I_j^{-1}\widetilde B(I_jx,\delta)\setminus F$ and connect $q_0$
and $q$ by a smooth oriented curve $\gamma\subset B(x,\delta)$ that
does not intersect $e$ and is transverse to $F$. By the properties of
$\delta$ and the above `leaves-and-enters' argument, we arrive at a
contradiction.

\medskip

{\bf Second step.} We are going to use Condition (3) in order to
`integrate' the above tessellation of a topological neighbourhood of
$P$ into a tessellation of a metric neighbourhood of $P$.

\smallskip

Fix some $\vartheta<d/2$, where $d$ is provided by IV, and fix some
$\varepsilon>0$ such that
$\varepsilon<\frac12\min\limits_{s\cap
s'\ne\varnothing}\varepsilon(s,s',\vartheta/2)$
and $\varepsilon<\vartheta/2$, where $\varepsilon(s,s',\vartheta/2)$ is
given by Condition (3).

Given an edge $e$, we put
$N_{e,r}:=\bigcup\limits_j\big(I_j^{-1},\widetilde N(I_je,r)\big)$,
where the $I_j$'s correspond to the geometric cycle including $e$ as in
3.1. For $s\in S$, define
$$N_{s,r}:=(1,\widetilde N(s,r)\big)\bigcup\big(I_{\overline
s},\widetilde N(\overline s,r)\big),\qquad W_s:=\pi
N_{s,\varepsilon}\bigcup_{e\in E_s}\pi N_{e,\vartheta}.$$

$\bullet$ Let us show that
$\varphi:W_s\to N(s,\varepsilon)\bigcup\limits_{e\in
E_s}N(e,\vartheta)$
is a bijection.

Choose any $e\in E_s$. As above, define
$F:=\bigcup\limits_jI_j^{-1}s_{j+1}\cap N(e,\vartheta)$, where the
$I_j$'s correspond to the geometric cycle including $e$ as in 3.1, and
pick a point $x\in e$. Using the tessellation of a small open ball $B$
centred at $x$, we can choose $q_0\in B$ living in exactly one of the
$I_j^{-1}\widetilde N(I_je,\vartheta)$. Clearly,
$F\subset\varphi\pi N_{e,\vartheta}$. It follows from the description
3.2 of formal neighbours that
$\varphi:\pi N_{e,\vartheta}\to N(e,\vartheta)$ is injective when
restricted to $F$. Let $q\in N(e,\vartheta)\setminus F$. As above,
connecting $q_0$ and $q$ by a smooth oriented curve
$\gamma\subset N(e,\vartheta)$ that does not intersect $e$ and is
transverse to $F$, we can see that $\gamma$ intersects only the
prescribed faces because $\vartheta<d$. We conclude that
$\varphi:\pi N_{e,\vartheta}\to N(e,\vartheta)$ is surjective and
injective. Since $\vartheta<d/2$, the $N(e,\vartheta)$ are disjoint.
Therefore,
$\varphi:\bigcup\limits_{e\in E_s}\pi
N_{e,\vartheta}\to\bigcup\limits_{e\in E_s}N(e,\vartheta)$
is a bijection.

It is easy to see that
$$s\setminus\bigcup_{e\in
E_s}N(e,\vartheta)\subset\varphi\Big(W_s\setminus\bigcup_{e\in E_s}\pi
N_{e,\vartheta}\Big)\subset N(s,\varepsilon)\setminus\bigcup_{e\in
E_s}N(e,\vartheta).$$
The description 3.2 of formal neighbours implies that
$\varphi:W_s\setminus\bigcup\limits_{e\in E_s}\pi N_{e,\vartheta}\to
N(s,\varepsilon)\setminus\bigcup\limits_{e\in E_s}N(e,\vartheta)$
is injective when restricted to
$s\setminus\bigcup\limits_{e\in E_s}N(e,\vartheta)$. Pick a point
$q\in N(s,\varepsilon)\setminus\bigcup\limits_{e\in E_s}N(e,\vartheta)$
such that $q\notin s$. There exist $x\in s$ and an oriented smooth
curve $\gamma\subset N(s,\varepsilon)$ of length
$\ell(\gamma)<\varepsilon$ that connects $x$ and $q$.

We claim that $\gamma$ can intersect $\partial P$ and
$\partial I_{\overline s}P$ only along $s\setminus\partial s$. Indeed,
$\gamma$ cannot intersect the faces of $P$ or of $I_{\overline s}P$
that are disjoint from $s$ because $\varepsilon<d$. Let $s'$ be a face
of $P$ or $I_{\overline s}P$ that intersects $\gamma$ and such that
$s\cap s'\ne\varnothing$. By Condition (3) and the choice of
$\varepsilon$, we have
$$s'\cap\gamma\subset s'\cap
N(s,\varepsilon)\subset\bigcup\limits_{s\diamond e\diamond
s'}N(e,\vartheta/2)\subset\bigcup\limits_{e\in E_s}N(e,\vartheta/2),$$
which implies $q\in\bigcup\limits_{e\in E_s}N(e,\vartheta)$ because
$\varepsilon<\vartheta/2$. A contradiction. The inequality
$\varepsilon<\vartheta/2$ implies that $\gamma$ does not intersect
$\partial s$.

We can assume that $\gamma$ is transverse to $s$. Considering the
tessellation of a small ball centred at $x$ introduced earlier, we see
that $\gamma$ first enters the interior of $P$ or of
$I_{\overline s}P$. When $\gamma$ intersects $s\setminus\partial s$, it
leaves $P$ and enters $I_{\overline s}P$ or {\it vice-versa.} As above,
$\varphi:W_s\setminus\bigcup\limits_{e\in E_s}\pi N_{e,\vartheta}\to
N(s,\varepsilon)\setminus\bigcup\limits_{e\in E_s}N(e,\vartheta)$
is surjective and injective.

\medskip

$\bullet$ Finally, let us show that the open metric neighbourhood
$N:=\overset{\,_\circ}\to P\bigcup\limits_{s\in
S}N(s,\varepsilon)\bigcup\limits_{e\in E_s}N(e,\vartheta)$
of $P$ is tessellated. (Note that $N(P,\varepsilon)\subset N$.) Define
$f(x)=\vartheta$ if $x\in e$ for some edge $e$ of $P$ and
$f(x)=\varepsilon$ if $x\in s\setminus\partial s$ for some $s\in S$. If
$x\in\overset{\,_\circ}\to P$, we take an arbitrary $f(x)>0$ such that
$B\big(x,f(x)\big)\subset\overset{\,_\circ}\to P$. It is immediate that
$\varphi W_{P,f}\subset N$ and that
$W_{P,f}=[\overset{\,_\circ}\to P]\bigcup\limits_{s\in S}W_s$. Hence,
$\varphi:W_{P,f}\to N$ is surjective. If $\varphi w=\varphi w'$, where
$w\in[\overset{\,_\circ}\to P]$ and $w'\in W_s$, then $w=[1,x]$,
$x\in\overset{\,_\circ}\to P$, and
$x\in N(s,\varepsilon)\bigcup\limits_{e\in E_s}N(e,\vartheta)$,
implying $w\in W_s$. If $x:=\varphi w=\varphi w'$, where $w\in W_s$ and
$w'\in W_{s'}$, then $s\ne s'$ and we have two cases:
$s\cap s'=\varnothing$ and $s\cap s'\ne\varnothing$. The first case is
impossible because
$N(s,\varepsilon)\bigcup\limits_{e\in E_s}N(e,\vartheta)\subset
N(s,\vartheta)$,
$N(s',\varepsilon)\bigcup\limits_{e\in E_{s'}}N(e,\vartheta)\subset
N(s',\vartheta)$,
and $N(s,\vartheta)\cap N(s',\vartheta)=\varnothing$ due to
$\vartheta<d/2$.

In the second case, suppose that $x\in N(e_0,\vartheta)$ for some
$e_0\in E_s$. Then $w\in\pi N_{e_0,\vartheta}$ because the bijection
$\varphi:\pi N_{e_0,\vartheta}\to N(e_0,\vartheta)$ is a restriction of
$\varphi:W_s\to N(s,\varepsilon)\bigcup \limits_{e\in
E_s}N(e,\vartheta)$
which is already known to be a bijection. Using
$\varepsilon<\vartheta<d/2$ and IV, we can see that the inclusion
$x\in N(s',\varepsilon)\bigcup\limits_{e\in E_{s}}N(e,\vartheta)$
implies $e_0\in E_{s'}$. So, $w'\in\pi N_{e_0,\vartheta}$ and $w=w'$.
The same arguments work if $x\in N(e_0,\vartheta)$ for some
$e_0\in E_{s'}$.

Therefore, we can assume that
$x\in N(s,\varepsilon)\cap N(s',\varepsilon)$ and
$x\notin\bigcup\limits_{s\diamond e\diamond s'}N(e,\vartheta)$. We can
find some
$p\in s'\cap B(x,\varepsilon)\subset s'\cap N(s,2\varepsilon)$. It
follows from $p\in B(x,\varepsilon)$,
$x\notin\bigcup\limits_{s\diamond e\diamond s'}N(e,\vartheta)$, and
$\varepsilon<\vartheta/2$ that
$p\notin\bigcup\limits_{s\diamond e\diamond s'}N(e,\vartheta/2)$. This
contradicts the choice of $\varepsilon$ and $\vartheta$ and completes
the proof
$_\blacksquare$

\bigskip

\centerline{\bf4.~Final Remarks}

\medskip

{\bf4.1.}~Probably, the first version of PPT where the restriction of
the combinatorial cycle isometry to an edge is not supposed to be the
identity can be found in [Kui, Subsection 3.1, p.~60], although in an
implicit form and in the specific case of real hyperbolic $4$-space.
(It is related to Remark 3.4.) Our~version seems to be the first
dealing with an angle condition in the situation of nonconstant angle
along a common edge of two faces.

\smallskip

{\bf4.2.}~Condition (3) in Theorem 3.5 is not trivial to check in the
case of nonconstant curvature. In~complex hyperbolic geometry, even in
so simple a case as that of bisectors intersecting transversally at a
common slice, the proof of this condition requires some analytic effort
[AGG, Lemma 2.2.3].

\smallskip

{\bf4.3.}~An important generalization of Theorem 3.5 would be of course
a version of PPT for polyhedra admitting faces of codimension $>2$.
Elaborating explicit conditions that express a good behaviour of faces
of all codimensions seems to be the most difficult task here. Indeed,
for simplicity, let us assume $M$ to be $3$-dimensional. Take a vertex
$p$ of $P$ and a small sphere $S$ centred at $p$. We have an
infinitesimal tessellation around generic points in edges, which
provides an infinitesimal tessellation of $S$ around its intersections
with the edges containing $p$. Due to the good behaviour of faces and
edges, we obtain a tessellation of $S$. While shrinking the radius of
$S$, the topological picture of this tessellation remains the same. In
this way, we visualize a tessellation of the $3$-ball bounded by $S$ as
being a cone over the tessellation of $S$.

In the particular case of a compact polyhedron, the conditions
expressing a good behaviour of faces must be drastically simplified.
For instance, the tessellation of a topological neighbourhood of the
polyhedron already implies in this case Tessellation Condition 2.1. We
suggest the following formulation: Let $M$ be a Riemannian manifold and
let $P\subset M$ be a simple compact PL-polyhedron equipped with a
face-pairing providing a family of geometric cycles of edges that
contain every edge of $P$. If Conditions~(1) and (2) of Theorem 3.5
hold, then Tessellation Condition 2.1 is satisfied.

This is just a rough outline of a possible proof; getting the general
version in question may require some serious effort. We thank Misha
Kapovich for pointing out the reference [Ale] (see also the proof of
[Ale, Theorem 2]) where similar ideas are applied to compact polyhedra
with totally geodesic faces in constant curvature spaces.

\newpage

\centerline{\bf5.~References}

\medskip

[AGG] S.~Anan$'$in, C.~H.~Grossi, N.~Gusevskii, {\it Complex hyperbolic
structures on disc bundles over surfaces,} to appear in
Int.~Math.~Res.~Not., see also http://arxiv.org/abs/math/0511741.

[Ale] A.~D.~Alexandrov, {\it On tiling a space with polyhedra,} in
A.~D.~Alexandrov Selected Works I, Classics of Soviet Mathematics
Vol.~4, Gordon and Breach Publishers, 1996. x+322 pp.

[Bea] A.~F.~Beardon, {\it The Geometry of Discrete Groups,} GTM
{\bf91}, Springer-Verlag, New York, 1983. xii+337 pp.

[EPe] D.~B.~A.~Epstein, C.~Petronio, {\it An exposition of Poincar\'e's
polyhedron theorem,} Enseign.~Math. {\bf40} (1994), 113--170.

[Kui] N.~H.~Kuiper, {\it Hyperbolic $4$-manifolds and tessellations,}
Inst.~Hautes \'Etudes Sci.~Publ.~Math. {\bf68} (1988), 47--76.

\enddocument